\documentclass[fleqn]{mat01}
\usepackage{times,mathtimy,amssymb,boxit}
\begin{document}
\setcounter{page}{179}
\firstpage{179}

\font\xx=msam5 at 9pt
\def\ab{\mbox{\xx{\char'03}}}

\font\sa=tibi at 10.4pt
\def\d{\hbox{d}}

\def\rem{\trivlist\item[\hskip\labelsep{{\it Remarks.}}]}
\def\thoe{\trivlist\item[\hskip\labelsep{{\bf Theorem.}}]}

\title{Beurling algebra analogues of the classical theorems\\ of Wiener
and L\'{e}vy on absolutely convergent\\ Fourier series}

\markboth{S J Bhatt and H V Dedania}{Analogues of Wiener's and
L\'{e}vy's theorems}

\author{S J BHATT and H V DEDANIA}

\address{Department of Mathematics, Sardar Patel University, Vallabh
Vidyanagar 388 120,  India\\
\noindent E-mail: subhashbhaib@yahoo.co.in; hvdedania@yahoo.com}

\volume{113}

\mon{May}

\parts{2}

\Date{MS received 20 February 2002; revised 26 September 2002}

\begin{abstract}
Let $f$ be a continuous function on the unit circle $\Gamma$, whose
Fourier series is $\omega$-absolutely convergent for some weight
$\omega$ on the set of integers $\mathcal{Z}$. If $f$ is nowhere
vanishing on $\Gamma$, then there exists a weight $\nu$ on $\mathcal{Z}$
such that $1/f$ had $\nu$-absolutely convergent Fourier series. This
includes Wiener's classical theorem. As a corollary, it follows that if
$\varphi$ is holomorphic on a neighbourhood of the range of $f$, then
there exists a weight $\chi$ on $\mathcal{Z}$ such that
\hbox{$\varphi\circ f$} has $\chi$-absolutely convergent Fourier series. This is a weighted
analogue of L\'{e}vy's generalization of Wiener's theorem. In the
theorems, $\nu$ and $\chi$ are non-constant if and only if $\omega$ is
non-constant. In general, the results fail if $\nu$ or $\chi$ is
required to be the same weight $\omega$.
\end{abstract}

\keyword{Fourier series; Wiener's theorem; L\'{e}vy's theorem; Beurling
algebra; commutative Banach algebra.}

\maketitle

\noindent Let $C(\Gamma)$ be the set of all continuous functions on the unit
circle $\Gamma$ in the complex plane $\mathcal{C}$. Let $f \in
C(\Gamma)$ such that the Fourier series
\begin{equation*}
f \sim \sum\limits_{n \in \mathcal{Z}} \widehat{f}(n) \
\hbox{e}^{int}\hbox{, where}\ \widehat{f} (n) = \frac{1}{2\pi i} \int_{-
\pi}^{\pi} f(\hbox{e}^{it}) \hbox{e}^{-int} \hbox{d} t\quad (n \in
\mathcal{Z}),
\end{equation*}
is absolutely convergent. If $f(z) \not=0$ for all $z \in \Gamma$, then
the Fourier series of $1/f$ is also absolutely convergent. This is a
classic Wiener's theorem (\cite{ERE}, \S11.4.17, p. 33), a transparent proof of
which by Gelfand (e.g. \cite{GRS}, p. 33) is often cited as the first success of
the theory of Banach algebras. L\'{e}vy's generalization of Wiener's
theorem states that if $\varphi$ is holomorphic on a neighbourhood of
the range of $f$, then $\varphi \circ f$ also has absolutely convergent
Fourier series (\cite{ERE}, \S11.4.17, p. 33). We aim to discuss Beurling
algebra analogues of these.

A {\it weight} on $\mathcal{Z}$ is a map $\omega : \mathcal{Z}
\longrightarrow [1,\infty)$ satisfying $\omega (m + n) \leq \omega (m)
\omega (n)$ for all $m, n \in \mathcal{Z}$. Let $\rho (1, \omega) =
\inf \{\omega (n)^{1/n} : n \geq 1\}$ and $\rho (2, \omega) = \sup
\{\omega(n)^{1/n}: n\leq - 1\}$. Then by (\cite{GRS}, p. 118), $0 < \rho (2,
\omega) \leq 1\leq \rho (1, \omega) < \infty$. A series $\sum_{n \in
\mathcal{Z}} \lambda_{n}$ is $\omega$-absolutely convergent if $\sum_{n
\in \mathcal{Z}} \big\vert\lambda_{n}\big\vert \omega(n) < \infty$. A
function $ f \in C(\Gamma)$ has $\omega$-absolutely convergent Fourier
series ($\omega$-ACFS) if its Fourier series is $\omega$-absolutely
convergent.

\begin{thoe}
{\it Let $\omega$ be a weight on $\mathcal{Z}$. Let $f\in C(\Gamma)${\rm ,} which
has $\omega$-ACFS.}
\end{thoe}\vspace{.4pc}

{\it
\noindent {\rm (I)}\ If $f(z) \not= 0$ for all $z \in \Gamma${\rm ,} then there exists a
weight $\nu$ on $\mathcal{Z}$ such that{\rm :}\vspace{.2pc}

\noindent {\rm (a)}\  $1/f $ has $\nu$-ACFS{\rm ;}

\noindent {\rm (b)}\  $\nu$ is non-constant if and only if $\omega$ is non-constant{\rm ;}

\noindent {\rm (c)}\  $\nu(n) \leq \omega(n)$ for all $n \in \mathcal{Z}$.\vspace{.2pc}

\noindent {\rm (II)}\  Let $\varphi$ be a function holomorphic on a neighbourhood of the
range of $f$. Then there exists a weight $\chi$ on $\mathcal{Z}$ such
that{\rm :}\vspace{.2pc}

\noindent {\rm (a)}\  $\varphi \circ f$ has $\chi$-ACFS{\rm ;}

\noindent {\rm (b)}\  $\chi$ is non-constant if and only if $\omega$ is non-constant{\rm ;}

\noindent {\rm (c)}\  $\chi (n) \leq \omega (n)$ for all $n \in \mathcal{Z}$.\vspace{.5pc}
}

The present note contributes to a programme suggested some thirty years
ago by Edward (\cite{ERE}, Ex. 11.15, p. 41). In the efforts made so far in this
programme, conditions on a given weight $\omega$ (e.g., the
Beurling--Domar condition; $\sum \frac{\log \omega(n)}{1 + n^{2}} < \infty$ (\cite{RS},
p. 185)) are sought, which ensure that $g$ (which is either $1/f$ or
$\varphi \circ f$ whatever the case may be) has $\omega$-ACFS. Contrary
to this, given an arbitrary weight $\omega$, we search for another
weight $\eta$ that ensure that $g$ has $\eta$-ACFS. We shall derive
(II) as a corollary of (I).

\begin{proof}
Let $\ell^{1} (\mathcal{Z}, \omega): = \big\{\lambda = (\lambda_{n}):
\big|\lambda \big|_{\omega} := \sum_{n\in \mathcal{Z}} \big|\lambda_{n}
\big| \omega(n) < \infty \big\} $, the Beurling algebra. It is a
convolution Banach algebra with norm $|\cdot|_{\omega}$. Let $A(\omega)
= \{ g \in C(\Gamma): \widehat{g} \in \ell^{1} (\mathcal{Z}, \omega)\}$,
the weighted Wiener algebra. It is a unital Banach algebra with the
pointwise operations and the norm being $\big\| g\big\|_{\omega} =
\big|\widehat{g}\big|_{\omega}$. Then $g\in C(\Gamma)$ has $\omega$-ACFS
if and only if $g \in A(\omega)$ and if and only if $\widehat{g} \in
\ell^{1} (\mathcal{Z}, \omega)$. Hence the Gelfand space $\Delta (A
(\omega))$ of $A(\omega)$ is identified with the closed annulus
$\Gamma(\omega) = \{z \in \mathcal{C} :\rho (2,\omega) \leq \big|z\big|
\leq \rho (1, \omega)\}$ via the map $z \in \Gamma(\omega)
\longmapsto \varphi_{z} \in \Delta (A(\omega))$, where $\varphi_{z} (g)
= \sum_{n \in \mathcal{Z}} \widehat{g} (n) z^{n} (g \in A(\omega))$.
Thus each function $g$ in $A(\omega)$ extends uniquely as an element
(denoted by $g$ itself) in $B(\omega)$ consisting of all continuous
functions on $\Gamma(\omega)$ which are analytic in its interior.\vspace{.8pc}

\noindent (I) Let $f \in C(\Gamma)$ have $\omega$-ACFS. Notice that $\Gamma
\subseteq \Gamma (\omega)$. Let $z \in \Gamma$. Since $f(z) \not= 0$,
there exists a neighbourhood $N(z)$ of $z$ in $\Gamma(\omega)$ such that
$\varphi_{w} (f) = f(w) \not=0$ for all $w \in N(z)$. We can assume that
$N(z) = \{w \in \mathcal{C}: \big|w - z\big| < r_{z}\} \cap
\Gamma(\omega)$ for some $r_{z} > 0$. By the compactness, there exist
$z_{1}, \ldots\,\!,z_{m}$ in $\Gamma$, arrange in such a way
that $\arg z_{i} < \arg z_{i + 1} \big(1 \leq i \leq m - 1\big)$, such
that $\Gamma \subseteq U_{1}^{m} N(z_{i}) \subseteq \Gamma(\omega)$. Now
we define positive numbers $r_{1}$ and $r_{2}$ as follows:\vspace{.5pc}

\noindent (i) If $\rho (2,\omega) = 1 = \rho (1, \omega)$, then take $r_{2} = 1
= r_{1}$.

\noindent (ii) If $\rho(2, \omega) = 1 < \rho (1, \omega)$, take $r_{2} = 1$; and
for $0 < \varepsilon < 1 - (1/\min \{s_{1},\ldots\!\,,s_{m}\})$, take
$r_{1} = (1 - \varepsilon) \min \{s_{1},\ldots\!\,,s_{m}\} > 1$, where
$s_{i} = \max \big\{\big| z\big| : z \in N(z_{i}) \cap N(z_{i + 1})\big\} (1
\leq i \leq m)$ and $z_{m + 1} = z_{1}$.

\noindent (iii) If $\rho (2, \omega) < 1 = \rho (1, \omega)$, take $r_{1} = 1$;
and for $0 < \varepsilon < (1/\max \{s_{1},\ldots\!\,,s_{m}\}) - 1$, take
$r_{2} = (1 + \varepsilon) \max \{s_{1},\ldots\!\,,s_{m}\} < 1$, where
$s_{i} = \min \big\{\big|z\big| : z \in N (z_{i}) \cap N(z_{i + 1})\big\} (1
\leq i \leq m)$ and $z_{m + 1} = z_{1}$.

\noindent (iv) If $\rho (2, \omega) < 1 < \rho (1, \omega)$, then take $r_{1}$
and $r_{2}$ as in (ii) and (iii) respectively.\vspace{.5pc}

Thus in any case, $\rho (2, \omega) \leq r_{2} \leq 1 \leq r_{1} \leq
\rho (1, \omega)$. Define $\nu : \mathcal{Z} \rightarrow [1, \infty)$ as
follows: If $\rho (2, \omega) = \rho (1, \omega)$, then take $\nu =
\omega$; otherwise define
\begin{equation*}
\nu(n) = \begin{cases}
r_{1}^{n} &\hbox{if} \ n \geq 0\\[.5pc]
r_{2}^{n} &\hbox{if} \ n \leq 0
\end{cases}.
\end{equation*}

It is clear that $\nu$ is non-constant if and only if $\omega$ is
non-constant. Then the following holds:\vspace{.3pc}

\noindent (1) $\nu$ is a weight on $\mathcal{Z}, \rho (2, \nu) = r_{2}$ and
$\rho (1, \nu) = r_{1}$;

\noindent (2) $\Gamma (\nu) \subseteq \Gamma (\omega)$;

\noindent (3) $f(z) \not= 0$ for all $z \in \Gamma(\nu)$;

\noindent (4) $1 \leq \nu (n) \leq \omega (n)$ for all $n \in \mathcal{Z}$.\vspace{.6pc}

Then by (4) above, $A(\omega)\subseteq A(\nu)$, and so $f\in A(\nu)$. Since $f(z)\neq 0$ for all $z$ in $\Gamma(\nu)=\Delta(A(\nu))$, it follows by the Gelfand theory that $1/f \in A(\nu)$, i.e. $1/f$ has $\nu$-ACFS.\vspace{.6pc}

\noindent (II) Let $K$ be the range of $f$. Let $\varphi$ be a function
holomorphic on a neighbourhood $U$ of $K$. Let $C$ be a closed
rectifiable Jordan contour in the open set $U$ containing $K$. Let $\mu
\in C$. Then $\mu \not\in K$ and $\mu 1 - f \in A(\omega)$. By part (I),
there exists a weight $\eta$ (which is non-constant if and only if
$\omega$ is non-constant) such that $\eta \leq \omega$ and the inverse
$(\mu1 - f)^{-1}$ of $(\mu 1- f)$ belongs to $A(\eta)$. Now take
$R_{\mu} = (\mu 1 - f)^{-1}$. Then its norm $\| R_{\mu} \|_{\eta}$ is
positive. Define $N(\mu) =\big\{ \lambda \in C : |\lambda - \mu| < \|
R_{\mu}\|_{\eta}^{-1}\}$. Then by the elementary Banach algebra argument,
it follows that for every $\lambda \in N(\mu),\lambda 1 - f = (\mu 1 -
f)\big\{ 1 + (\lambda  -\mu) R_{\mu}\big\}$ is invertible in $A(\eta)$.
Thus $N(\mu)$ is a neighbourhood of $\mu$ in $C$ such that for all
$\lambda \in N(\mu), \lambda {1}  - f$ is invertible in $A(\eta)$.

Now by the compactness of $C$, there exist finitely many
$\mu_{1},\ldots\,\!, \mu_{n}$ in $C$ and weights
$\eta_{1},\ldots\eta_{n}$ such that $C \subseteq \cup_{1}^{n}
N(\mu_{i})$, and for any $\lambda \in C$, the inverse of $\lambda 1 - f$
belongs to $A(\eta_{i})$ for some $i$. Now define
\begin{equation*}
r_{2} = \max \big\{\rho (2, \eta_{i}) : 1 \leq i \leq n\big\} \
\hbox{and} \ r_{1} = \min \big\{ \rho (1, \eta_{i}) : 1 \leq i \leq n
\big\}
\end{equation*}
so that $r_{2} \leq 1 \leq r_{1}$. If $\rho (2, \omega) = 1 = \rho (1,
\omega)$, then by Part~I, each $\eta_{i} = \omega$. If $\rho (2, \omega)
= 1 < \rho (1, \omega)$, then $\rho (2, \eta_{i}) = 1 < \rho (1,
\eta_{i})$ for each $i$, and so $r_{2} = 1 < r_{1}$. Similarly, the
cases $\rho(2, \omega) < 1 = \rho (1, \omega)$ and $\rho(2,\omega) < 1
<\rho (1, \omega)$ can be discussed. Now if $\rho (2, \omega) = 1 = \rho
(1, \omega)$, then take $\chi = \omega( = \eta_{i})$; otherwise define
$\chi : \mathcal{Z} \longrightarrow [1, \infty)$ as
\begin{equation*}
\chi(n) = \begin{cases}
r_{1}^{n} &\hbox{if}\ n \geq 0\\[.5pc]
r_{2}^{n} &\hbox{if} \ n \leq 0
\end{cases}.
\end{equation*}
It is clear that $\chi$ is non-constant if and ony if $\omega$ is non-constant. Then the following holds.
\begin{enumerate}
\renewcommand{\labelenumi}{(\arabic{enumi})}
\item $\chi$ is a weight on $\mathcal{Z}, \rho (2, \chi) = r_{2}$ and
$\rho (1,\chi) = r_{1}$;

\item $\rho (2,\omega) \leq \rho (2, \eta_{i}) \leq \rho (2, \chi) \leq
1 \leq \rho (1, \chi) \leq \rho (1, \eta_{i}) \leq \rho (1, \omega)$ for
all $i$;

\item $1 \leq \chi \leq \eta_{i} \leq \omega$ on $\mathcal{Z}$ and hence
$A(\omega) \subseteq A(\eta_{i}) \subseteq A(\chi)$ for all $i$;

\item For any $\lambda \in C$, the inverse of $\lambda 1 - f$ belongs to
$A(\chi)$.\vspace{-.55pc}
\end{enumerate}

Now the map $\lambda \in C \longrightarrow \varphi (\lambda)
R_{\lambda}$ is a continuous map from $C$ into the Banach algebra $\big(
A(\chi), \|\cdot\|_{\chi}\big)$, where $R_{\lambda}$ is the
inverse of $\lambda 1 - f$. Hence the integral $(1/2\pi i) \int_{C}
\varphi (\lambda) R_{\lambda} \hbox{d}\lambda$ is in $A (\chi)$ in the
sense of $\|\cdot\|_{\chi}$-convergence and $\varphi(f) =
(1/2\pi i) \int_{C} \varphi (\lambda) R_{\lambda} \hbox{d}\lambda$, where
$\varphi (f)$ is defined by the functional calculus in $C(\Gamma)$. Thus
$\varphi (f)$ has $\chi$-ACFS. It follows that $\varphi (f)
\big(\hbox{e}^{i\theta}\big) = \big(\varphi \circ f\big)
\big(\hbox{e}^{i\theta}\big)$ for all $\hbox{e}^{i\theta}\in \Gamma$.\hfill $\Box$
\end{proof}

\begin{rem}$\left.\right.$\vspace{.3pc}

\noindent (1) Let $\omega$ be any weight on $\mathcal{Z}$ such that $\rho
(2,\omega) \not= \rho (1,\omega)$. Then $\Gamma$ is properly contained
in $\Gamma(\omega)$. Let $f\in C(\Gamma)$ have $\omega$-ACFS such that
$f(z) \not= 0$ for all $z \in \Gamma$, and $f(z_{0}) = 0$ for some $z_{0}
\in \Gamma(\omega)$. Then the function $f$ is clearly not invertible in
$A(\omega)$, i.e., $1/f$ cannot have $\omega$-ACFS. For example,
define $\omega (n) = e^{|n|} \ (n \in \mathcal{Z})$ and let $f(z) =
z_{0} - z (z \in \mathcal{C})$, where $1 < |z_{0}| < e$. Then $f$ has
$\omega$-ACFS, $\rho (1, \omega) = e, \rho (2,\omega) = 1/e$ and $1/f$
does not have $\omega$-ACFS.

\noindent (2) Let $\omega$ be a weight on $\mathcal{Z}$ such that $\rho (2,
\omega) = 1 = \rho (1, \omega)$. Then it follows from the proof that for
any $f\in C(\Gamma)$ having $\omega$-ACFS and satisfying $f(z) \not= 0$
for all $z \in \Gamma$, the $1/f$ has also $\omega$-ACFS. Examples of such
weights include:

\noindent (i) $\omega_{\alpha} (n) = \big(1 + |n|\big)^{\alpha}$, where $0 < \alpha <
\infty$;

\noindent (ii) $\omega (n) = 1 + \log \big(1 + |n|\big)$;

\noindent (iii) $\omega (n) = \big(1 + |n|\big)^{\sqrt{1 + |n|}}$.

\noindent (3) Let $f \in C(\Gamma)$ such that $f$ have $\omega$-ACFS for
every weight $\omega$ on $\mathcal{Z}$. Suppose $f(z) \not= 0$ for all
$z \in \Gamma$. One would be tempted to know whether $1/f$ has $\omega$-ACFS for every $\omega$. The answer is `no'. For example, take $f(z) =
2z + z^{2}$, a trigonometric polynomial. Then the Fourier series of
$1/f$ is
\begin{equation*}
\left(\frac{1}{f}\right) (z) = \frac{1}{2z} \sum\limits_{0}^{\infty} (-
1)^{k} \left(\frac{z}{2}\right)^{k}
\end{equation*}
which fails to have $\omega$-ACFS for the weight $\omega(n) = 2^{|n| +
2} (n \in \mathcal{Z})$.
\end{rem}

\section*{Acknowledgement}

The authors are thankful to the referee for a careful reading of the
manuscript.


\begin{thebibliography}{99}
\bibitem{ERE} Edwards R-E, Fourier Series, (New York: Holt, Rinehart and
Winston Inc.) (1967) vol. II

\bibitem{GRS} Gelfand I, Raikov D and Shilov G, Commutative normed rings
(New York: Chelse Publication Company) (1964)

\bibitem{RS} Reiter H and Stegeman J D, Classical harmonic analysis and
locally compact abelian groups (Oxford: Clarendon Press) (2000)
\end{thebibliography}
\end{document}